\documentclass[10pt,twocolumn,twoside]{IEEEtran}
\newcommand{\beq}{\begin{equation}}
\newcommand{\eeq}{\end{equation}}
\newcommand{\bsq}{\begin{subequations}}
\newcommand{\esq}{\end{subequations}}
\newcommand{\bq}{\begin{eqnarray}}
\newcommand{\eq}{\end{eqnarray}}
\newcommand{\bqn}{\begin{eqnarray*}}
\newcommand{\eqn}{\end{eqnarray*}}
\usepackage{cite}
\usepackage{url}
\usepackage{empheq}
\usepackage{float}
\usepackage{dsfont}
\usepackage[hidelinks]{hyperref}

\usepackage{threeparttable}
\usepackage{algorithm}
\usepackage{algorithmic}
\usepackage{braket}

\usepackage[table,xcdraw]{xcolor}

\usepackage{arydshln}
\usepackage{nicematrix}

\usepackage{amsmath,amssymb,amsfonts}
\allowdisplaybreaks[4]
\usepackage{graphicx}
\usepackage{textcomp}
\usepackage{amsmath}
\usepackage{enumitem}
\usepackage{epstopdf}
\usepackage{array}
\usepackage{booktabs}
\usepackage{subfigure}
\usepackage{multirow}
\usepackage{soul, color}
\usepackage[utf8]{inputenc}
\usepackage{cases}
\usepackage{amsthm}  

\newtheorem{theorem}{Theorem}

\newtheorem{definition}{Definition}
\newtheorem{remark}{Remark}

\newtheorem{assumption}{Assumption}

\soulregister\cite7
\soulregister\citep7
\soulregister\citet7
\soulregister\ref7
\soulregister\pageref7

\usepackage{bbding}

\usepackage{nomencl}
\makenomenclature
\usepackage{etoolbox}
\renewcommand\nomgroup[1]{%
  \item[\bfseries
  \ifstrequal{#1}{A}{Acronyms}{%
  \ifstrequal{#1}{S}{Symbols}{%
  \ifstrequal{#1}{U}{Units}{}}}%
]}

\makeatletter
\newcommand{\Statex}{\item[]}                                           \makeatother

\NiceMatrixOptions{code-for-first-row = \scriptstyle,code-for-first-col = \scriptstyle }

\newif\ifincludeSecondAppendix
\includeSecondAppendixtrue

\begin{document}

\title{Quantum-Enabled Probabilistic Optimal Power Flow with Built-in Differential Privacy}

\author{Yuji Cao,~\IEEEmembership{Graduate Student Member,~IEEE},
    Tongxin Li,
    Yue Chen,~\IEEEmembership{Senior Member,~IEEE}
    %
    \thanks{Y. Cao and Y. Chen are with the Department of Mechanical and Automation Engineering, The Chinese University of Hong Kong, Hong Kong, China. (email: \{yjcao, yuechen\}@mae.cuhk.edu.hk)}
    \thanks{T. Li is with the School of Data Science, The Chinese University of Hong Kong, Shenzhen, China. (email: litongxin@cuhk.edu.cn).}
}

\markboth{IEEE Transactions on Smart Grid}%
{Cao et al.: Quantum-Enabled Probabilistic Optimal Power Flow with Built-in Differential Privacy}

\maketitle

\begin{abstract}
Quantum computing has been regarded as a promising approach to accelerate power system optimization. However, challenges such as limited qubits and inherent noise hinder their widespread adoption in power systems. In this paper, we propose a qubit-efficient framework for solving a crucial power system optimization problem, the probabilistic optimal power flow (POPF). We demonstrate that quantum noise, traditionally viewed as a drawback, can in fact be leveraged to provide a built-in differential privacy (DP) guarantee. Specifically, we first linearize POPF into a multi-parametric linear program (MP-LP) with renewable uncertainties being the parameters. This decomposes the parameter space into critical regions with precomputed solution maps. Second, a variational quantum circuit (VQC) classifies the critical region based on each uncertainty realization and then recovers the final solution. In this way, the required qubits scale with the uncertain parameters instead of the network size, with only 5 qubits versus 600+ for direct quantum OPF in a 69-bus system. Moreover, we prove the depolarizing noise of VQC provides DP guarantees and characterize the privacy-cost tradeoff. Case studies validate the proposed VQC achieves 2.1$\times$ smaller privacy budgets compared to its classical counterpart. At matched privacy levels, the VQC also maintains lower infeasibility and prediction error.

\end{abstract}











\section{Introduction}

\IEEEPARstart{T}{he} proliferation of renewable energy sources is driving a global transition towards decarbonized power grids \cite{zakariaUncertaintyModelsStochastic2020}. While these sources can lower carbon emissions, they also introduce significant volatility and uncertainty. To address this challenge, probabilistic optimal power flow (POPF) is now a mainstream analytical framework. Its objective is to assess the statistical characteristics of critical operating quantities under uncertainty. However, POPF requires evaluating thousands to millions of scenarios and leads to a prohibitive computational burden for conventional methods~\cite{liNonparametricProbabilisticOptimal2022}.

Quantum computing (QC) has emerged as a potential solution to such computationally intractable problems in power systems~\cite{morstynOpportunitiesQuantumComputing2024}. Unlike classical methods, QC exploits quantum-mechanical principles such as superposition and entanglement for computation. For problems with certain structures, QC is expected to provide exponential speedup over classical computers~\cite{lee2023evaluating}. These prospects have motivated growing interest in applying QC to power systems optimization~\cite{han2025quantum}. However, hardware challenges remain in current noisy intermediate-scale quantum (NISQ) devices.

A major hardware challenge is the limited number of quantum bits, i.e., qubits. Typical OPF problems involve hundreds of continuous decision variables. Solving these on quantum hardware requires discretizing each continuous variable into a binary representation, where the number of qubits per variable determines the precision. For hundreds of OPF variables, the required number of qubits far exceeds current hardware capacity~\cite{pareekDemystifyingQuantumPower2024}. To address this problem, one line of work~\cite{morstyn2022annealing} focuses on combinatorial OPF with discrete decision variables only and solves it via quantum annealing. Another work~\cite{nikmehr2022quantum} decomposes the problem so that quantum hardware solves only the discrete subproblems, while continuous variables are handled by classical solvers. The quantum speedup is thus limited to the combinatorial component of the problem.

The other major hardware challenge is the inevitable quantum noise~\cite{wangNoiseinducedBarrenPlateaus2021,caiQuantumErrorMitigation2023}. Indeed, imperfections in quantum hardware introduce errors that accumulate during computation and degrade the final results. To suppress the noise, techniques such as quantum error mitigation and correction have been developed~\cite{temmeErrorMitigationShortDepth2017,googlequantumaiQuantumErrorCorrection2024}. Despite these advances, full error correction still requires substantial additional qubits and further limits the feasibility for near-term OPF applications. Meanwhile, near-term mitigation techniques offer only limited protection as circuit depth grows. As a result, non-negligible hardware noise remains unavoidable on near-term devices and must be accounted for when designing quantum optimization methods.

Rather than viewing the noise solely as a liability, recent work in information theory and statistics has shown that the intrinsic randomness of quantum noise can itself protect privacy~\cite{zhouDifferentialPrivacyQuantum2017,hircheQuantumDifferentialPrivacy2022}. This observation is particularly relevant to the optimal power flow problem, which requires sensitive operational data such as customer loads and generator costs as input. OPF solutions, in turn, can be exploited to infer such private data~\cite{dvorkin2021dpopf}. To provide a rigorous privacy guarantee, differential privacy (DP) mechanisms have been widely applied to OPF problems, typically through artificial noise injection~\cite{ryu2021privacy,zhou2019differential}. By contrast, NISQ devices are inherently affected by hardware noise. This raises a natural question: does such built-in hardware noise already contribute to a meaningful DP guarantee for OPF, and if so, how strong is that contribution? Without answering this question, adding extra artificial noise on top of hardware noise is unnecessarily conservative and further degrades OPF solution quality.


To address the above challenges, this paper proposes a quantum-enabled framework for probabilistic OPF. We reduce the qubit requirements by decomposing the parameter space into critical regions and using a variational quantum circuit (VQC) to identify the active region. Moreover, instead of viewing inherent quantum noise as a liability, we use it as a natural random resource, prove the DP guarantee and characterize the privacy-cost tradeoff. Our main contributions are two-fold:
\begin{enumerate}
    \item \textbf{A qubit-efficient quantum-enabled POPF framework:} We propose a quantum-enabled POPF framework built on multi-parametric linear programming (MP-LP) critical-region decomposition. The framework precomputes affine solution maps for each critical region and uses a VQC to classify the active region during operation. This reduces the qubit requirement from encoding all decision and slack variables to classifying the active critical region. On a 69-bus system, the proposed framework requires only 5 qubits, compared with $596\sim1406$ qubits required for direct encoding of OPF formulations at different precisions, with an estimated per-scenario speedup of over $9,000\times$.

    \item \textbf{Privacy guarantees from inherent quantum noise:} 
    We develop a depolarizing noise model for the proposed quantum-enabled POPF framework and prove that it provides $(\varepsilon, 0)$-differential privacy. We further bound the expected OPF cost increase under randomized dispatch and derive an explicit privacy-cost tradeoff. Case studies on a 69-bus system show that the VQC achieves $2.1\times$ tighter privacy bounds than a classical multi-layer perceptron (MLP). At matched privacy levels ($\varepsilon \approx 3$), the proposed quantum-enabled framework maintains 0.21\% cost gap and 4.6\% infeasibility, compared to 16.07\% and 66.1\% for the MLP with Gaussian noise.
\end{enumerate}

The remainder of this paper is organized as follows. Section~\ref{sec:problem} formulates POPF via multi-parametric programming. Section~\ref{sec:qnn} presents the quantum classifier. Section~\ref{sec:privacy-analysis} analyzes privacy guarantees and the privacy-cost tradeoff. Case studies are conducted in Section~\ref{sec:case-studies} and conclusions are drawn in Section~\ref{sec:conclusion}.
\section{Problem Formulation}
\label{sec:problem}

In this section, we first formulate the probabilistic optimal power flow (POPF) problem to be addressed. Then, we introduce a multi-parametric linear programming (MP-LP) approach that partitions the parameter space into critical regions, within which the optimal solution can be expressed as an affine function of the parameters. This reformulation transforms the POPF solving process into identifying the critical region corresponding to the given parameters, which significantly reduces the number of qubits required for computation.

\subsection{Probabilistic OPF and Linearization}

The POPF problem seeks to calculate the expected minimum operational cost of power systems over uncertainty, such as renewable generation and load demand. For each uncertainty scenario, denoted by $\theta$, a deterministic OPF is solved. The POPF problem is formulated as:
\begin{subequations}
\label{eq:prob_opf}
\begin{align}
    \mathbb{E}_{\theta} \Big[ \min_{u}\ & f(u, \theta) \Big], \label{eq:prob_opf-obj} \\
    \text{s.t.} \quad & h(u, \theta) = 0, \label{eq:prob_opf-eq} \\
    & g(u, \theta) \leq 0, \label{eq:prob_opf-ineq}
\end{align}
\end{subequations}
where $u$ denotes the decision variables (e.g., thermal generator setpoints, voltage magnitudes), $\theta$ is a vector of uncertain parameters, and $f(\cdot)$ is the operational cost. $h(\cdot)=0$ and $g(\cdot)\leq 0$ encode the equality and inequality operational constraints, respectively. Evaluating the expectation in~\eqref{eq:prob_opf} typically requires Monte Carlo sampling, where the inner optimization problem is solved repeatedly for a large number of uncertainty scenarios.

The inner optimization for each scenario is an OPF problem. Consider a network with a bus set $\mathcal{N}$, a generator set $\mathcal{G}$, and a line set $\mathcal{L}$. Let $\mathcal{G}_i \subseteq \mathcal{G}$ denote generators at bus $i$. The standard alternating current (AC) OPF problem that minimizes the total generation cost subject to power balance and operational limits is given by
\begin{subequations}
\label{eq:ACOPF}
\begin{align}
    \min_{P_g, Q_g, V, \delta} ~ & \sum_{g \in \mathcal{G}} C_g(P_g),
    \label{eq:ac-opf-obj} \\
\text{s.t.} ~ & V_i \sum_{j \in \mathcal{N}} V_j \Bigl( G_{ij} \cos \delta_{ij}  + B_{ij} \sin \delta_{ij} \Bigr) \nonumber \\
& = \sum_{g \in \mathcal{G}_i} P_g - P_{d,i} ,\forall i \in \mathcal{N}, \label{eq:ac-opf-p}\\
& V_i \sum_{j \in \mathcal{N}} V_j \Bigl( G_{ij} \sin \delta_{ij} - B_{ij} \cos \delta_{ij} \Bigr)\nonumber\\
& = \sum_{g \in \mathcal{G}_i} Q_g - Q_{d,i}, \forall i \in \mathcal{N},\label{eq:ac-opf-q}  \\
& \underline{P}_g \le P_g \le \overline{P}_g, \; \underline{Q}_g \le Q_g \le \overline{Q}_g, \; \forall g \in \mathcal{G}, \label{eq:ac-opf-glim} \\
& \underline{V}_i \le V_i \le \overline{V}_i, \; \forall i \in \mathcal{N}, \label{eq:ac-opf-vlim}
\end{align}
\end{subequations}
where $P_g$ and $Q_g$ are the active and reactive power outputs of generator $g$, $P_{d,i}$ and $Q_{d,i}$ are the active and reactive demands at bus $i$, and $V_i$ and $\delta_i$ are the voltage magnitude and angle with $\delta_{ij} := \delta_i - \delta_j$. The terms $G_{ij}$ and $B_{ij}$ are elements of the bus admittance matrix $Y = G + jB$, and $C_g(\cdot)$ is the generation cost function. 
Constraints~\eqref{eq:ac-opf-p}--\eqref{eq:ac-opf-q} enforce power balance at each bus, while~\eqref{eq:ac-opf-glim}--\eqref{eq:ac-opf-vlim} impose operational limits.

The power flow equations \eqref{eq:ac-opf-p} and \eqref{eq:ac-opf-q} make the AC-OPF problem nonconvex. Various linearization methods have been developed to obtain tractable approximations, including direct current (DC) power flow models~\cite{stott2009dc}, LinDistFlow for radial networks \cite{baran1989network}, and second-order cone programming relaxation \cite{low2014convex1} plus polyhedral approximations \cite{chen2018energy}. When such linearization is applied, the AC-OPF model \eqref{eq:ACOPF} reduces to a linear programming (LP) problem.

Although the linearized OPF is an LP, evaluating the expectation in~\eqref{eq:prob_opf} still requires one LP per uncertainty scenario. For large-scale Monte Carlo simulations with thousands of scenarios, this repeated computation remains expensive. Quantum computing has the great potential to accelerate such workloads, but current NISQ devices provide only tens to hundreds of qubits. The full OPF decision vector may contain hundreds of variables, far exceeding this capacity. To facilitate the implementation with limited qubits, in the following, we propose a multi-parametric programming-based reformulation method, which reduces the OPF problem to a classification problem over a low-dimensional parameter space.


\subsection{Multi-parametric Linear Programming Reformulation}

When the inner OPF problem is linearized, treating the uncertainties as parameters yields a multi-parametric linear programming (MP-LP). Let $x \in \mathbb{R}^n$ denote the decision variables and $\theta \in \mathbb{R}^m$ the uncertain parameters. The problem can be expressed in a compact parametric form:
\begin{equation}
    \label{eq:MP-LP}
    \begin{aligned}
        x^*(\theta) = \arg \min_{x} \quad & c^\top x  ,                                                 \\
        \text{s.t.} \quad                                                        & Wx \le S + T\theta,
    \end{aligned}
\end{equation}
where $c \in \mathbb{R}^n$ is the cost coefficient vector, $W \in \mathbb{R}^{q \times n}$ and $S \in \mathbb{R}^q$ encode the constraint coefficients, and $T \in \mathbb{R}^{q \times m}$ determines how the uncertain parameters $\theta$ enter the constraint right-hand side. The parameter $\theta$ ranges over a bounded polytope $\Theta \subset \mathbb{R}^m$ determined by feasible load and/or generation bounds. Equality constraints (e.g., power balance) are handled by elimination or by expressing them as paired inequalities.

This MP-LP reformulation recasts the computational task. A key result from MP-LP theory is that the parameter space $\Theta$ can be partitioned into a finite number (say $K$) of convex polyhedral sets, known as critical regions (CRs). Within each critical region $k$, denoted by $\text{CR}_k$, the set of active constraints at the optimal solution remains invariant. This invariance implies that the optimal solution vector $x^*(\theta)$ becomes a direct affine function of the parameter $\theta$:
\begin{equation}
    \label{eq:critical-region}
    x^*(\theta) = F_k \theta + f_k, \quad \forall \theta \in \text{CR}_k.
\end{equation}
To obtain $F_k$ and $f_k$, consider the set of active constraints within $\text{CR}_k$. Let $\mathcal{A}_k$ denote the index set of constraints that are binding at the optimal solution. These constraints satisfy
\begin{equation}
    \label{eq:active-constraints}
    W_{\mathcal{A}_k} x^*(\theta) = S_{\mathcal{A}_k} + T_{\mathcal{A}_k}\theta,
\end{equation}
where $W_{\mathcal{A}_k}$, $S_{\mathcal{A}_k}$, and $T_{\mathcal{A}_k}$ are the rows of $W$, $S$, and $T$ indexed by $\mathcal{A}_k$. Assuming the active constraint matrix has full column rank and, in the non-degenerate case, consists of $n$ linearly independent constraints, $W_{\mathcal{A}_k} \in \mathbb{R}^{n \times n}$ is nonsingular. Therefore,
\begin{equation}
    \label{eq:Fr}
    x^*(\theta)
    = W_{\mathcal{A}_k}^{-1}S_{\mathcal{A}_k}
    + W_{\mathcal{A}_k}^{-1}T_{\mathcal{A}_k}\theta.
\end{equation}
Applying \eqref{eq:critical-region} to \eqref{eq:Fr}, it follows that
\[
F_k = W_{\mathcal{A}_k}^{-1}T_{\mathcal{A}_k}, \qquad
f_k = W_{\mathcal{A}_k}^{-1}S_{\mathcal{A}_k}.
\]
The inverse formula above applies in the non-degenerate case, where the active matrix contains $n$ linearly independent binding constraints. In a degenerate region, more than $n$ constraints may be binding, and multiple primal or dual optimal solutions may make the construction of critical regions from optimal basis invariance difficult. Following \cite{chen2022flexibility}, in this paper we first recover the critical region from a dual or value-function representation and then recover the corresponding primal solution. Therefore, the piecewise-affine representation in \eqref{eq:critical-region} is preserved.

The critical region decomposition transforms the solution of the OPF problem into a more manageable two-stage process. First, a multi-class classification is performed: for a given uncertainty scenario $\theta$, one must identify the critical region $\text{CR}_k$ it belongs to. Since each critical region corresponds to a unique active constraint set $\mathcal{A}_k$, the binding constraints are directly determined once the region index $k$ is known. Second, once the region is identified, the full, high-dimensional optimal solution $x^*(\theta)$ is recovered with negligible computation via the simple affine solution map in \eqref{eq:critical-region}. This transformation shifts the computational burden from high-dimensional regression to classification. Rather than regressing all $n$ OPF decision variables, the problem reduces to identifying the region index $k$, from which the full solution is recovered analytically via \eqref{eq:critical-region}. This classification is well-suited to quantum computing because quantum feature maps embed the high-dimensional uncertainty parameter $\theta$ into an exponentially large Hilbert space \cite{havlivcek2019supervised}, providing richer representational power for separating critical regions than classical methods of comparable complexity.

\section{Hybrid Quantum-Classical Classifier for Probabilistic OPF}
\label{sec:qnn}

This section presents the hybrid quantum-classical classifier for fast critical region identification. In the classifier, a variational quantum circuit architecture extracts features from the encoded uncertainty parameters. A classical network then maps these features to region scores, and a temperature-controlled softmax produces region probabilities, enabling tunable tradeoffs between classification confidence and output randomness.

\subsection{Variational Quantum Feature Extractor}

\label{sec:vqc}
To identify the critical region that contains a given uncertainty scenario $\theta$, we employ a variational quantum classifier. Its fundamental component is the variational quantum circuit (VQC), which uses a parameterized quantum circuit as a trainable feature map to learn complex representations in the high-dimensional Hilbert space of the quantum system~\cite{cerezo2021variational}. The VQC maps the uncertainty $\theta$ to a set of distinguishing features to separate the $K$ critical regions.

\begin{figure}[bhtp]
    \centering
    \includegraphics[width=\linewidth]{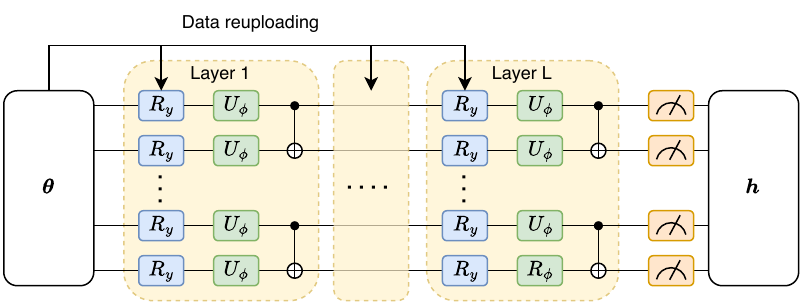}
    \caption{Variational quantum feature extractor with data reuploading technique. Each layer re-encodes input $\theta$ via $R_y$ gates, followed by trainable $U_\phi$ gates and CNOT entanglement. Measurements on all qubits yield the feature vector $h$.}
    \label{fig:vqc-architecture}
\end{figure}

As illustrated in Fig.~\ref{fig:vqc-architecture}, the VQC transforms the classical input $\theta \in \mathbb{R}^m$ into the feature vector $h$ through a layered quantum circuit followed by measurement. The circuit acts on $n_q$ qubits initialized in the ground state $\ket{0}^{\otimes n_q}$, producing an output state $\ket{\psi_{\text{out}}}$ from which classical information is extracted. The ket notation $\ket{\cdot}$ denotes a quantum state vector in a $2^{n_q}$-dimensional space. Each feature $h_i$ is the Pauli-$Z$ expectation value on qubit $i$, a real number in $[-1, 1]$ indicating whether qubit $i$ is more likely measured in state $\ket{0}$ (yielding $+1$) or $\ket{1}$ (yielding $-1$):
\begin{equation}
  h_i = \braket{\psi_{\text{out}} | \sigma_z^{(i)} | \psi_{\text{out}}}, \quad i = 1, \dots, n_q.
\end{equation}

A conventional VQC architecture encodes the classical input data only once initially, then applies trainable gates. However, such single-encoding circuits have limited expressivity: their outputs are restricted to Fourier series with frequencies determined solely by the encoding gates~\cite{schuld2021effect}. This limits the circuit's ability to learn complex, nonlinear decision boundaries required for accurate region classification.

To address this limitation, we adopt the \emph{data reuploading} technique~\cite{perez2020data}, which interleaves data encoding with trainable gates across multiple layers. By re-encoding $\theta$ at each layer, the circuit can represent richer function classes with higher-frequency Fourier components~\cite{schuld2021effect}. This design is particularly suitable for near-term quantum devices, as it achieves greater expressivity without requiring additional qubits.

Specifically, the circuit in Fig. \ref{fig:vqc-architecture} comprises $L$ repeated layers. Each layer $\ell$ consists of three operations applied sequentially:
\begin{enumerate}
  \item \textbf{Data encoding:} The input $\theta$ is encoded via single-qubit $R_y$ rotations. For qubit $i$, the gate is defined as:
      \begin{equation}
          R_y(\theta_i) = \exp\left(-i \frac{\theta_i}{2} \sigma_y\right) = \begin{pmatrix} \cos\frac{\theta_i}{2} & -\sin\frac{\theta_i}{2} \\ \sin\frac{\theta_i}{2} & \cos\frac{\theta_i}{2} \end{pmatrix},
      \end{equation}
      where $\theta_i$ is determined by a corresponding component of $\theta$ (with appropriate preprocessing when $m \neq n_q$).
      \item \textbf{Trainable rotations:} Parameterized single-qubit gates $R_y(\phi_i^{(\ell)})$ are applied, where $\phi_i^{(\ell)}$ are learnable parameters optimized during training. This enables the circuit to learn task-specific transformations.
      \item \textbf{Entanglement:} CNOT gates arranged in a ladder topology create correlations between neighboring qubits:
      \begin{equation}
          U_{\text{ent}} = \prod_{i=1}^{n_q-1} \text{CNOT}_{i, i+1},
      \end{equation}
      where $\text{CNOT}_{i,j}$ flips qubit $j$ conditioned on qubit $i$. This entanglement allows the circuit to capture multi-feature interactions that cannot be represented by product states.
  \end{enumerate}

  The full circuit is given by the composition of all layers:
  \begin{equation}
      U(\phi_q, \theta) = \prod_{\ell=1}^{L} \left[ U_{\text{ent}} \cdot U_\phi^{(\ell)} \cdot U_E(\theta) \right],
  \end{equation}
  where $U_E(\theta) = \bigotimes_{i=1}^{n_q} R_y(\theta_i)$ denotes the data-encoding layer, and $U_\phi^{(\ell)} = \bigotimes_{i=1}^{n_q} R_y(\phi_i^{(\ell)})$ the trainable layer. The output state is:
  \begin{equation}
      \ket{\psi_{\text{out}}} = U(\phi_q, \theta) \ket{0}^{\otimes n_q}.
      \label{eq:psi_out}
  \end{equation}

Finally, classical information can be extracted from the final quantum state $\ket{\psi_{\text{out}}}$, by measuring the Pauli-Z expectation value on each qubit. For qubit $i$, the feature is:
  \begin{equation}
      h_i = \braket{\psi_{\text{out}} | \sigma_z^{(i)} | \psi_{\text{out}}} = p_i^{(0)} - p_i^{(1)},
      \label{eq:measurement}
  \end{equation}
where $p_i^{(0)}$ and $p_i^{(1)}$ are the probabilities of measuring qubit $i$ in state $\ket{0}$ and $\ket{1}$, respectively. In practice, these expectation values are estimated by repeated circuit executions (shots) and averaging the measurement outcomes.

This quantum-extracted vector $h$ encodes a rich feature representation of the input $\theta$, which includes high-order correlations learned by the VQC and thus provides a powerful basis for the subsequent network.

\subsection{Temperature-Controlled Region Classifier}
Building on representations learned from the VQC, a classical network maps the compact representation to region scores and learns weights that determine how features are combined to distinguish each region.  A softmax function with inverse temperature $\beta$ then converts these scores to selection probabilities. Once a region is selected, the full OPF solution is reconstructed via the precomputed affine solution map \eqref{eq:critical-region}.


The VQC feature vector $h$ is passed to a linear classifier that produces logits $s = W_{cl} h + b_{cl}$, where $W_{cl}$ and $b_{cl}$ are weights and biases. A softmax with inverse temperature $\beta > 0$ converts these logits to region selection probabilities:
  \begin{equation}
      \label{eq:softmax}
      p_k = \frac{\exp(\beta s_k)}{\sum_{j=1}^{K} \exp(\beta s_j)}.
  \end{equation}
The parameter $\beta$ controls the sharpness of the output distribution: large $\beta$ concentrates probability mass on the top-scoring region for accurate classification, while small $\beta$ flattens the distribution toward uniform. During training, quantum parameters $\phi_q$ and classical parameters $\phi_c = \{W_{cl}, b_{cl}\}$ are trained jointly by minimizing cross-entropy loss, with quantum gradients obtained via the parameter-shift rule~\cite{wierichs2022general} below:

\begin{equation}
  \frac{\partial \mathcal{L}}{\partial \phi} = \frac{1}{2} \left[ \mathcal{L}(\phi + \pi/2)
  - \mathcal{L}(\phi - \pi/2) \right].
\end{equation}
These gradients are then fed to a classical optimizer such as Adam to update all model parameters.



\begin{algorithm}[tbhp]
    \caption{A Hybrid Quantum-Classical POPF Algorithm}
    \label{alg:pipeline}
    \begin{algorithmic}[1]
        \STATE \textbf{Offline Phase: Training}
        \REQUIRE Training data $\{(\theta_i, k_i)\}_{i=1}^N$, affine solution maps $\{(F_k, f_k)\}_{k=1}^K$, Epochs $E$, batch size $B$, learning rate $\eta$
        \STATE Initialize $\phi_q$ and $\phi_c = \{W_{cl}, b_{cl}\}$
        \FOR{epoch $= 1$ to $E$}
            \FOR{each mini-batch $\mathcal{B}$}
                \STATE Forward:
                \Statex $U(\phi_q, \theta)\ket{0}^{\otimes n_q} \to h \to \text{softmax}(W_{cl}h + b_{cl})$
                \STATE Compute cross-entropy loss $\mathcal{L}$
                \STATE Gradients: backprop ($\phi_c$) + parameter-shift ($\phi_q$)
                \STATE Update: $\phi_q \leftarrow \phi_q - \eta \nabla_{\phi_q} \mathcal{L}$, \; $\phi_c \leftarrow \phi_c - \eta \nabla_{\phi_c} \mathcal{L}$
            \ENDFOR
        \ENDFOR
        \STATE Store trained parameters $(\phi_q^*, \phi_c^*)$
    \end{algorithmic}
    \hrule
    \begin{algorithmic}[1]
        \STATE \textbf{Online Phase: Privacy-Preserving Inference}
        \REQUIRE New scenario $\theta_{\text{new}}$, temperature $\beta$, noise level $\gamma$
        \STATE Process with noise:
        \Statex \hspace{1.5em} $\rho_{\text{out}} = \mathcal{E}_\gamma\bigl(U(\phi_q^*, \theta_{\text{new}})\ket{0}\bra{0}^{\otimes n_q}U^\dagger(\phi_q^*, \theta_{\text{new}})\bigr)$
        \STATE Measure: $h_j = \mathrm{Tr}(\sigma_z^{(j)} \rho_{\text{out}})$ for $j = 1, \ldots, n_q$
        \STATE Compute logits: $s = W_{cl}^* h + b_{cl}^*$
        \STATE Compute probabilities: $p_k = \exp(\beta s_k) / \sum_{j=1}^{K} \exp(\beta s_j)$
        \STATE Sample region: $\tilde{k} \sim \text{Categorical}(p)$
        \STATE Reconstruct solution: $x^* = F_{\tilde{k}} \theta_{\text{new}} + f_{\tilde{k}}$
        \ENSURE Region index $\tilde{k}$ (released), OPF solution $x^*$
    \end{algorithmic}
\end{algorithm}

The complete offline training and online inference procedures are detailed in Algorithm~\ref{alg:pipeline}. Once trained, the model offers a substantial acceleration for POPF analysis. For any new uncertainty scenario $\theta_{\text{new}}$, the hybrid classifier rapidly infers the corresponding critical region index $k$. Subsequently, the full, high-dimensional OPF solution vector $x^*(\theta_{\text{new}})$ is reconstructed with negligible computation on a classical computer via the affine transformation in~\eqref{eq:critical-region}. This paradigm addresses the dimensionality challenge of near-term quantum devices by shifting the computational task from high-dimensional regression to low-dimensional classification.

However, deploying this quantum classifier on near-term hardware raises two practical concerns. First, noise in current quantum hardware distorts the measured features $h$ and degrades classification accuracy. Second, the released region index $\tilde{k}$ may leak sensitive information about the operating condition $\theta$, such as individual load patterns or renewable generation. A key observation is that these two concerns are inherently coupled. The same quantum noise that degrades classification accuracy also obscures the mapping from $\theta$ to $\tilde{k}$, which is a natural mechanism for privacy protection. In the next section, we formalize this observation and quantify the resulting privacy-accuracy tradeoff.

\section{Privacy Analysis}
\label{sec:privacy-analysis}

Hardware noise is an inherent feature of near-term quantum devices. This section shows that such noise provides differential privacy for the released region index. We first define the inherent quantum noise and formalize how depolarizing noise affects the output distribution, then prove the main privacy guarantee, and finally characterize the privacy-cost tradeoff.

\subsection{Inherent Quantum Noise}
\label{sec:privacy-model}

We assume an adversary observes the region index $\tilde{k}$ released by the online dispatch pipeline and attempts to infer sensitive information about the current operating point $\theta$ (e.g., customer loads or renewable generation). Our objective is to ensure that the distribution of $\tilde{k}$ is insensitive to small perturbations of $\theta$, formalized by the following adjacency relation.

\begin{definition}[Parameter Adjacency]
    \label{def:adjacency}
    For a fixed $\Delta_\theta > 0$, two parameter vectors $\theta, \theta' \in \mathbb{R}^m$ are adjacent (denoted by $\theta \sim \theta'$) if $\|\theta - \theta'\|_2 \le \Delta_\theta$.
\end{definition}

The radius $\Delta_\theta$ bounds the perturbation for which outputs should remain indistinguishable. In power systems, this corresponds to operationally similar scenarios, such as a 5\% fluctuation in renewable generation or load demand. Our goal is to provide privacy guarantees to ensure that an adversary cannot distinguish such neighboring conditions by observing the released region index.

In this paper, we argue that the inherent quantum noise can obscure the released region index and thus protect privacy. First,
we model quantum hardware errors using the depolarizing channel, a standard noise model that approximates the aggregate effect of gate errors in near-term quantum processors~\cite{urbanek2021mitigating}. For any quantum state $\rho$ on $n_q$ qubits, the depolarizing channel with noise level $\gamma \in (0,1)$ is defined as follows:
\begin{equation}
    \mathcal{E}_\gamma(\rho) = (1-\gamma)\rho + \gamma\frac{I}{D}.
\end{equation}
That is, with probability $1-\gamma$, the state is preserved; with probability $\gamma$, it is replaced by $I/D$, with $D = 2^{n_q}$ and $I \in \mathcal{R}^{D \times D}$ is the identity matrix. The state $I/D$ represents complete loss of information, in which all measurement outcomes are equally likely.

Let $\rho(\theta)$ denote the density matrix of the quantum state after encoding the input $\theta$ via the data-reuploading circuit. Under this noise model, the measured feature on qubit $j$ becomes
\begin{equation}
    h_j(\theta) = \mathrm{Tr}\left( \sigma_z^{(j)} \, \mathcal{E}_\gamma\left( U(\phi_q)\, \rho(\theta)\, U^\dagger(\phi_q) \right) \right).
\end{equation}
The noise level $\gamma$ controls the contraction strength. Larger $\gamma$ flattens the output distribution toward uniform, which strengthens privacy but reduces classification accuracy.

We make the following assumption:

\begin{assumption}[Lipschitz Encoding]
    \label{ass:lipschitz}
    There exists $L_{\mathrm{enc}} \ge 0$ such that for all $\theta, \theta'$,
    \begin{equation}
        D_{\mathrm{tr}}(\rho(\theta), \rho(\theta')) \le L_{\mathrm{enc}} \|\theta - \theta'\|_2,
    \end{equation}
    where $D_{\mathrm{tr}}(\rho, \rho') := \frac{1}{2}\|\rho - \rho'\|_1$ is the trace distance.
\end{assumption}

The trace distance $D_{\mathrm{tr}}(\rho, \rho')$ quantifies how distinguishable two quantum states are, ranging from zero (identical) to one (orthogonal). Assumption~\ref{ass:lipschitz} requires that nearby parameters produce similar quantum states, with $L_{\mathrm{enc}}$ bounding the encoding's sensitivity. This assumption is satisfied by common encoding schemes such as angle encoding~\cite{sim2019expressibility, berberich2024training}, where the rotation angles depend continuously on the input. It links parameter-space adjacency to quantum-state proximity.

\subsection{Differential Privacy Guarantee}
\label{sec:dp-guarantee}

In the following, we analyze how the inherent quantum noise can protect privacy. We first present the concept of differential privacy.

\begin{definition}[Differential Privacy~\cite{dworkRoth2014}]
    \label{def:dp}
    A randomized mechanism $\mathcal{M}: \Theta \to \mathcal{Y}$ satisfies $(\varepsilon, 0)$-differential privacy if for any adjacent $\theta \sim \theta'$ and any measurable set $S \subseteq \mathcal{Y}$,
    \begin{equation}
        \mathbb{P}[\mathcal{M}(\theta) \in S] \le e^\varepsilon \, \mathbb{P}[\mathcal{M}(\theta') \in S].
    \end{equation}
\end{definition}

Based on this definition, the following theorem quantifies the privacy budget $\varepsilon$ of Algorithm 1 in terms of the depolarizing noise level $\gamma$ and inverse temperature $\beta$.

\begin{theorem}[Region-ID Differential Privacy]
    \label{thm:dp}
    Under Assumption~\ref{ass:lipschitz}, the randomized mechanism in Algorithm~\ref{alg:pipeline} that maps $\theta$ to the released region index $\tilde{k}$ satisfies $(\varepsilon_{\mathrm{reg}}, 0)$-differential privacy over the output space 
    $\mathcal{Y} = \{1, \ldots, K\}$, where
    \begin{equation}
        \label{eq:eps-bound}
        \varepsilon_{\mathrm{reg}} \le 4\beta(1-\gamma)\, L_{\mathrm{enc}} \Delta_\theta \, \|W_{cl}\|_{\infty,1},
    \end{equation}
    with $\|W_{cl}\|_{\infty,1} := \max_{k} \|w_k\|_1$ denoting the maximum $\ell_1$-norm of the rows of the weight matrix $W_{cl}$ in the linear classifier.
\end{theorem}

The proof of Theorem~\ref{thm:dp} can be found in Appendix~\ref{appendix:thm-dp}. It shows that depolarizing noise uniformly contracts quantum expectation values, reducing the sensitivity of output probabilities to input perturbations.


\subsection{Privacy-Cost Tradeoff}
\label{sec:tradeoff}

The privacy guarantee in Theorem~\ref{thm:dp} comes at a cost: randomized region selection may result in an incorrect region, leading to a suboptimal dispatch. This section quantifies the expected cost increase and establishes an explicit tradeoff between privacy and operational efficiency.

Let $J(x; \theta)$ denote the operational cost when solution $x$ is applied under parameter $\theta$. For any region $k$, the affine solution from~\eqref{eq:critical-region} is $x_k(\theta) = F_k \theta + f_k$, and the classifier logit from~\eqref{eq:softmax} is $s_k(\theta)$. Let $k^*(\theta)$ denote the correct region containing $\theta$, with optimal solution $x^*(\theta) = F_{k^*} \theta + f_{k^*}$.

Define the cost gap for selecting region $k$ as $\Delta J_k(\theta) := J(x_k(\theta); \theta) - J(x^*(\theta); \theta)$. Under randomized dispatch, the realized cost gap is $\Delta J := \Delta J_{\tilde{k}}(\theta)$, where $\tilde{k}$ is the sampled region index. The score margin $m(\theta) := s_{k^*}(\theta) - \max_{k \ne k^*} s_k(\theta)$ measures the logit margin between the correct region and the highest-scoring competing region. A larger positive margin indicates clearer separation.



\begin{theorem}[Privacy-Cost Tradeoff]
    \label{thm:tradeoff}
    Let $\Delta J_{\max}(\theta) := \max_{k \ne k^*} \Delta J_k(\theta)$. Under softmax sampling with inverse temperature $\beta$,
    \begin{equation}
        \label{eq:tradeoff}
        \mathbb{E}[\Delta J \mid \theta] \le \underbrace{\Delta J_{\max}(\theta)}_{\text{worst-case cost gap}} \cdot \underbrace{(K-1) \exp(-\beta \, m(\theta))}_{\text{mis-selection probability bound}},
    \end{equation}
    where $\Delta J_{\max}(\theta)$ captures the maximum cost penalty from selecting a wrong region, while $(K-1) \exp(-\beta \, m(\theta))$ bounds the probability of such mis-selection via the softmax margin $m(\theta)$.
\end{theorem}

The proof of Theorem~\ref{thm:tradeoff} is provided in the Appendix~\ref{appendix:thm-tradeoff}. It proceeds in three parts: first establishing the cost decomposition, then bounding the mis-selection probability via the softmax margin, and finally combining them to derive the tradeoff.

\begin{remark}[Simplified Bound for Standard VQC]
    \label{rem:vqc-simplification}
    For standard variational quantum circuits using Pauli-$Z$ measurements and no classifier bias ($b_{cl} = 0$), the margin scales linearly with noise: $m(\theta) = (1-\gamma) \, m^{(0)}(\theta)$, where $m^{(0)}(\theta)$ is the noiseless margin. This yields a simplified bound:
    \begin{equation}
        \label{eq:tradeoff-vqc}
        \mathbb{E}[\Delta J \mid \theta] \le \Delta J_{\max}(\theta) \cdot (K-1) \exp\left( -\beta (1-\gamma) \, m^{(0)}(\theta) \right). \nonumber
    \end{equation}
\end{remark}

\begin{remark}[Tuning the Privacy-Cost Tradeoff]
    \label{rem:tuning}
    The privacy budget $\varepsilon_{\mathrm{reg}}$ and expected cost increase are both controlled by $\beta(1-\gamma)$. Increasing $\gamma$ (more noise) or decreasing $\beta$ (softer sampling) improves privacy but raises cost. Eliminating $\beta(1-\gamma)$ between Theorems~\ref{thm:dp} and~\ref{thm:tradeoff} yields an explicit relationship: for standard VQC (Remark~\ref{rem:vqc-simplification}),
    \begin{equation}
        \mathbb{E}[\Delta J \mid \theta] \le \Delta J_{\max}(\theta) \cdot (K-1) \exp\left( -\frac{m^{(0)}(\theta)}{4 L_{\text{enc}}\Delta_\theta \|W_{cl}\|_{\infty,1}} \, \varepsilon_{\mathrm{reg}} \right).  \nonumber
    \end{equation}
    Given a privacy budget $\varepsilon_{\mathrm{reg}}$, the cost bound decays exponentially with $\varepsilon_{\mathrm{reg}}$, and the decay rate improves with larger margin $m^{(0)}$.
\end{remark}

\begin{figure}[htbp]
    \centering
    \includegraphics[width=\columnwidth]{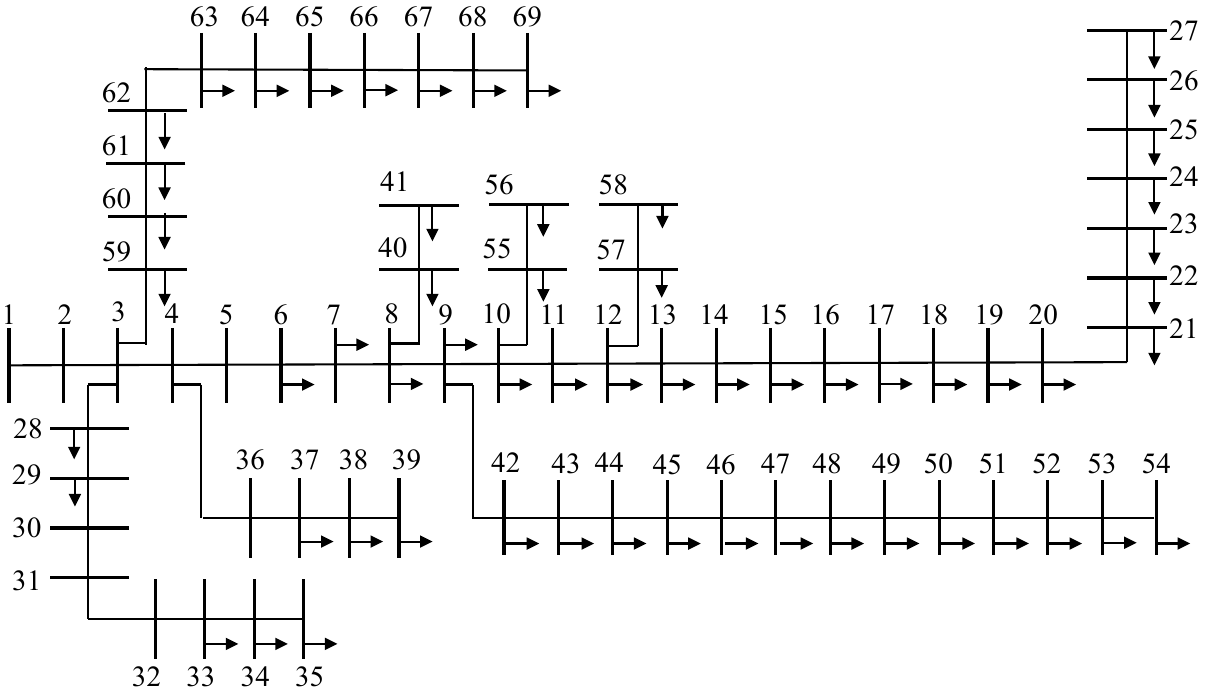}
    \caption{Topology of 69-bus distribution network.}
    \label{fig:topology}
\end{figure}

\begin{figure}[bthp]
    \centering
    \includegraphics[width=0.8\columnwidth]{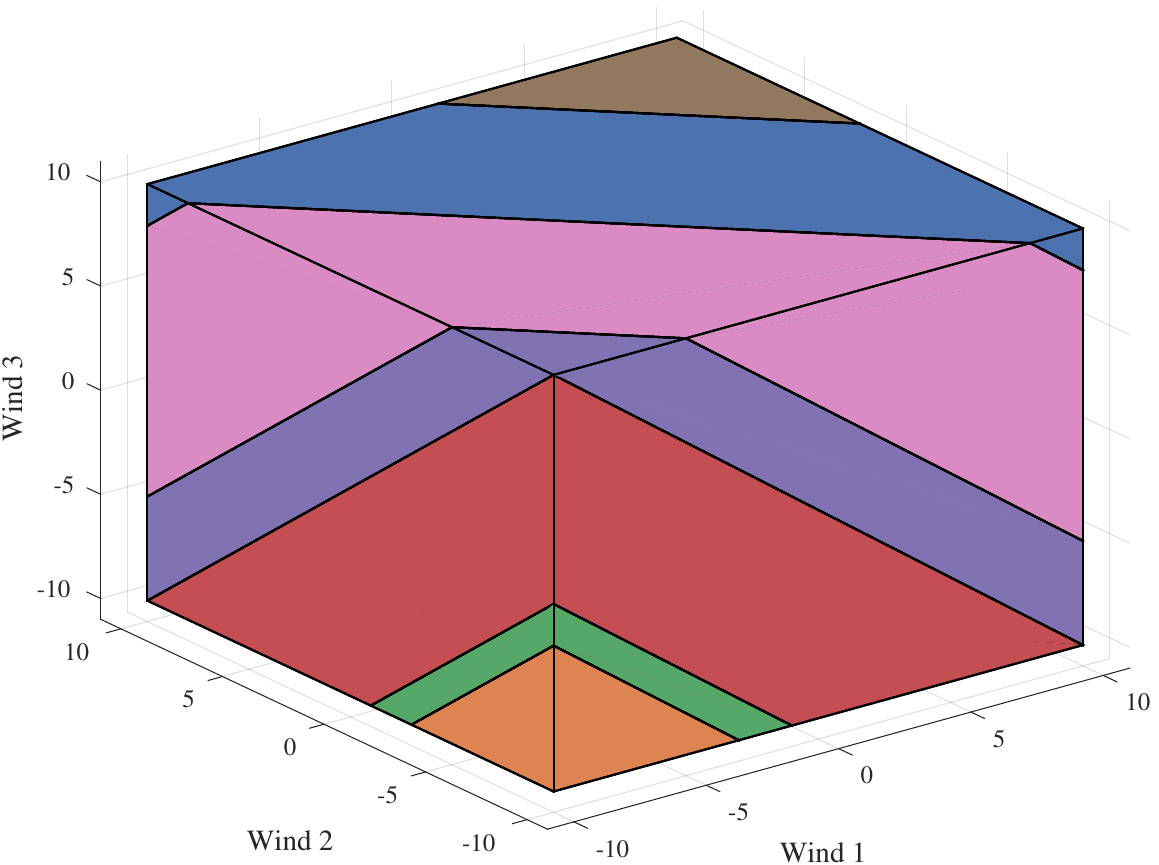}
    \caption{The critical regions of the multi-parametric LP problem for the 69-bus distribution network.}
    \label{fig:critical_region}
\end{figure}

\section{Case Study}
\label{sec:case-studies}


We test our method on a modified 69-bus distribution network~\cite{kadir2013optimal}. 
The topology is shown in Fig.~\ref{fig:topology}. Six elastic demands are located at nodes 12, 23, 32, 42, 53 and 62, respectively. Three renewable generators are connected to nodes 9, 30, and 60, respectively. Their real-time deviations fall within $[-10, 10]$ kW. Using these deviations as the uncertainty parameter $\theta = [\Delta w_1, \Delta w_2, \Delta w_3]^\top$ and applying the MP-LP method described in~\eqref{eq:MP-LP}, we obtain seven critical regions as shown in Fig.~\ref{fig:critical_region}. 3,000 data points are sampled from the uncertainty space and used to train the variational quantum classifier with uncertainty parameters normalized to $[-1, 1]$.

\begin{table}[bhtp]
    \centering
    \caption{Model Architecture and Performance Comparison}
    \label{tab:model_config}
    \begin{tabular}{lcc}
        \toprule
                           & VQC                & MLP                     \\
        \midrule
        Architecture       & 5 qubits, 6 layers & $3 \to 7 \to 7 \to 7$ \\
        Parameters         & 125                & 133                     \\
        Test Accuracy (\%) & 98.23              & 97.78                   \\
        \bottomrule
    \end{tabular}
\end{table}

A VQC and a multi-layer perceptron (MLP) are trained to predict the region index based on the uncertainty parameters. Table~\ref{tab:model_config} summarizes the model architectures and their performance. The VQC uses 5 qubits with 6 variational layers, each consisting of angle encoding gates followed by strongly entangling CNOT gates. The MLP uses 2 hidden layers with 7 neurons each, matching the VQC's parameter count. A linear layer without bias is appended to both models. Both models are trained for 30 epochs using Adam optimizer (learning rate 0.05, batch size 32) with cross-entropy loss under noise-free conditions. The VQC is implemented in PennyLane~\cite{bergholm2018pennylane} and the MLP in PyTorch~\cite{paszke2019pytorch}.

\subsection{Privacy Guarantees and Trade-offs}

\begin{figure}[hbtp]
    \centering
    \includegraphics[width=1.0\columnwidth]{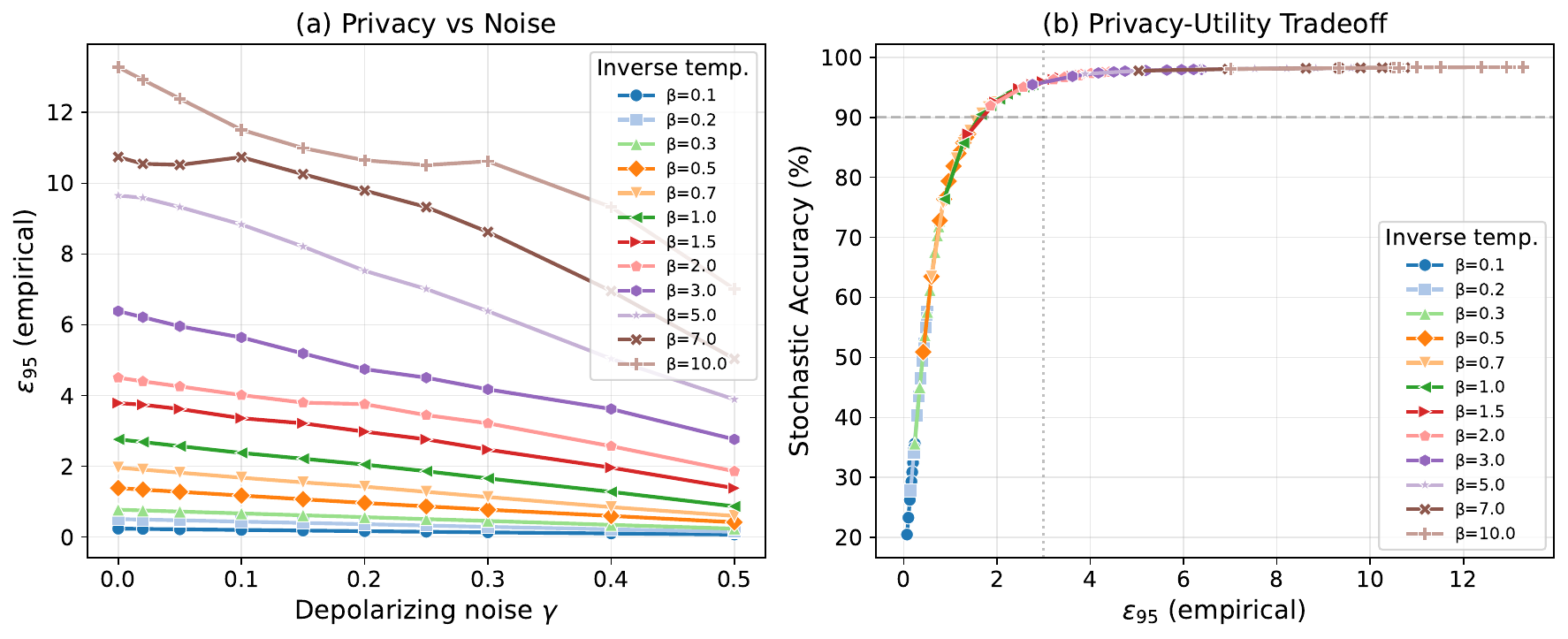}
    \caption{Privacy-utility frontier showing stochastic accuracy versus empirical $\varepsilon$ for different inverse temperatures $\beta$. Each curve traces the effect of varying noise $\gamma$ from 0 to 0.5.}
    \label{fig:privacy_utility}
\end{figure}

\begin{figure*}[hbtp]
    \centering
    \includegraphics[width=\linewidth]{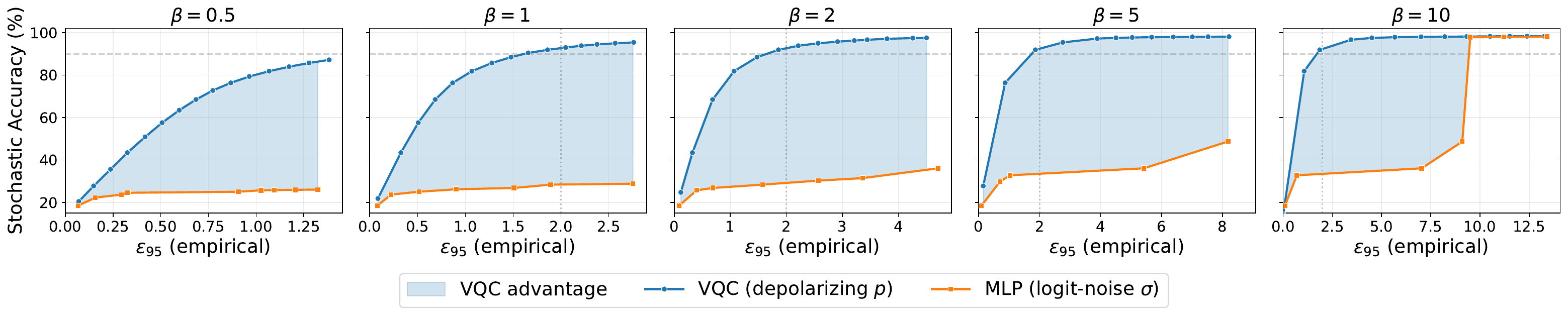}
    \caption{Privacy-utility frontier comparison between VQC (with depolarizing noise) and multi-layer perceptron (MLP, with Gaussian logit noise). At matched privacy levels, the VQC maintains substantially higher accuracy.}
    \label{fig:classical_comparison}
\end{figure*}

The randomized region-selection mechanism provides tunable privacy through two parameters: the depolarizing noise level $\gamma$ and the inverse temperature $\beta$. This subsection validates the differential privacy guarantee established in Theorem~\ref{thm:dp}, characterizes the resulting privacy-utility tradeoff, and compares the VQC against a classical MLP baseline.

To empirically assess the privacy budget, we construct 100 adjacent input pairs and measure how the output distribution changes. The adjacent parameters $\theta$ and $\theta'$ are bounded by $\Delta_\theta = 0.05$ with $\|\theta - \theta'\|_2 \le \Delta_\theta$ in the normalized feature space, corresponding to approximately 1~kW perturbation in load or renewable generation. For each pair, we compute the empirical privacy budget
\begin{equation}
    \varepsilon_{\mathrm{emp}}(\theta, \theta') = \max_{k} \left| \log \frac{P_k(\theta)}{P_k(\theta')} \right|,
\end{equation}
which measures the maximum log-ratio of output probabilities. We report $\varepsilon_{95}$, the 95th percentile over adjacent pairs, as an empirical privacy indicator that is robust to outliers. Note that $\varepsilon_{95}$ is a statistical summary of observed log-ratios, not a formal DP guarantee; the formal worst-case bound is $\varepsilon_{\mathrm{reg}}$ from Theorem~\ref{thm:dp}.


The left figure in Fig.~\ref{fig:privacy_utility} analyzes the effects of noise and temperature. We evaluate $\varepsilon_{95}$ for $\gamma \in \{0, 0.1, 0.2, 0.3, 0.4, 0.5\}$ and $\beta$ from 0.1 to 10. As predicted by Theorem~\ref{thm:dp}, increasing depolarizing noise $\gamma$ reduces the privacy budget $\varepsilon_{95}$ due to the $(1-\gamma)$ factor. Similarly, decreasing the inverse temperature $\beta$ flattens the softmax distribution, further reducing $\varepsilon_{95}$. The right figure reveals the resulting privacy-utility tradeoff by evaluating the accuracies for all $(\gamma, \beta)$ combinations. The result collapses into a single frontier. The frontier exhibits a steep transition region: accuracy degrades rapidly below $\varepsilon_{95} \approx 3$, where privacy constraints begin to dominate, but remains above 90\% for $\varepsilon_{95} \ge 3$. This characterization allows operators to analyze the tradeoff between privacy and accuracy and select an optimal parameter configuration based on their privacy requirements.

\begin{figure}[htbp]
    \centering
    \includegraphics[width=0.9\columnwidth]{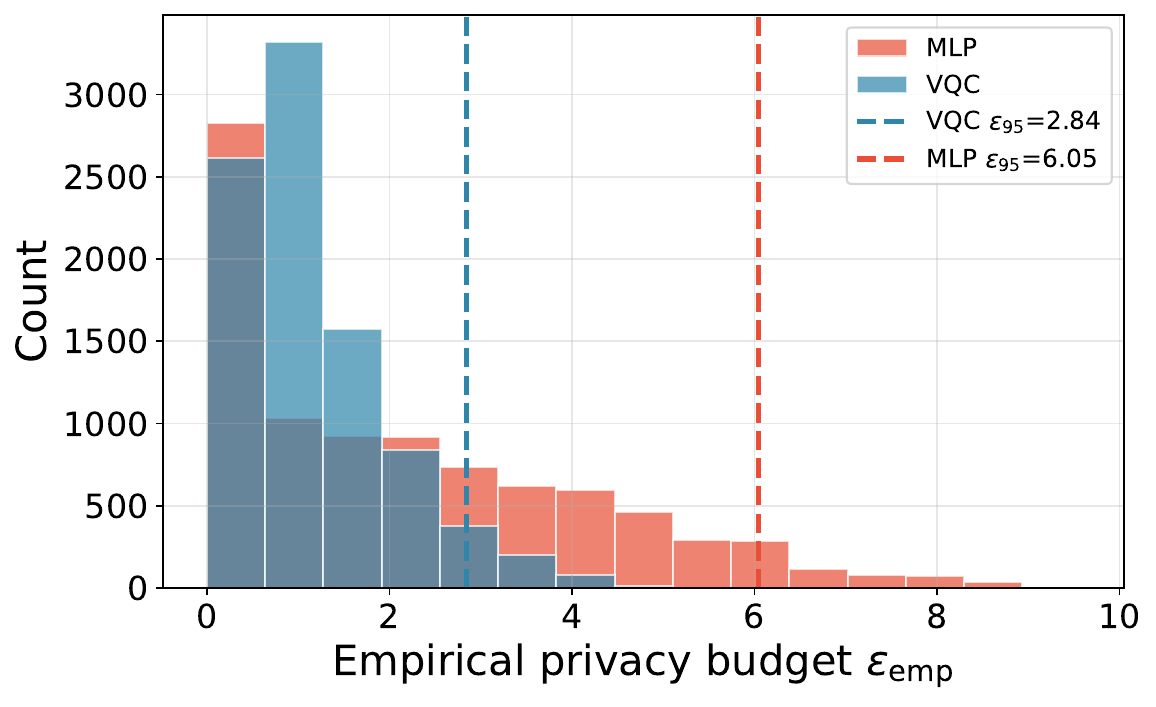}
    \caption{Distribution of empirical privacy budget $\varepsilon_{\mathrm{emp}}$ over adjacent pairs. The VQC (blue) has a shorter tail than the MLP (red), yielding $\varepsilon_{95} = 2.84$ vs. $6.05$.}
    \label{fig:logratio_distribution}
\end{figure}

Fig.~\ref{fig:classical_comparison} compares the privacy-utility frontiers of the VQC and a classical MLP baseline across different inverse temperatures. The MLP baseline shares the same classical head architecture as the VQC. To introduce privacy, independent Gaussian noise $\mathcal{N}(0, \sigma^2)$ is added to the MLP logits before softmax sampling. For each target privacy level, the noise scale $\sigma$ is tuned until the resulting $\varepsilon_{95}$ matches the target, using the same evaluation protocol (100 adjacent pairs, $\Delta_\theta = 0.05$) as for the VQC. In the figures, we can observe a gap (shaded area) between the accuracy of VQC and MLP, indicating that the VQC achieves higher accuracy than the MLP at matched $\varepsilon_{95}$. The advantage is more notable at low $\beta$: at $\beta = 1$ and $\varepsilon_{95} \approx 1$, the VQC achieves approximately 80\% accuracy while the MLP remains near 27\%, which is barely above random guessing for seven critical regions. Even at higher $\beta$ values where the MLP eventually reaches high accuracy, it requires significantly larger $\varepsilon_{95}$. This consistent gap demonstrates that depolarizing noise provides a more favorable privacy-utility tradeoff than additive Gaussian noise on classical logits.

To explain this advantage, we examine the baseline case without depolarizing noise. At $\gamma=0$ and $\beta=1$, the VQC achieves $\varepsilon_{95} = 2.84$ while the MLP yields $\varepsilon_{95} = 6.05$, a $2.1\times$ difference even without any privacy-enhancing mechanism. Fig.~\ref{fig:logratio_distribution} visualizes the distribution of empirical privacy budget over adjacent pairs: the VQC distribution concentrates at lower values with 95\% below 2.84, while the MLP tail extends beyond 6.

\begin{table}[hbtp]
  \centering
  \begin{threeparttable}
  \caption{Probability Distribution Sensitivity to Input Perturbation}
  \label{tab:confidence_comparison}
  \begin{tabular}{cccccc}
  \toprule
  & & \multicolumn{2}{c}{\textbf{Probability} $P_{k'}$} & & \\
  \cmidrule(lr){3-4}
  \multirow{-2}{*}{\textbf{Method}} & \multirow{-2}{*}{$P_{k^*}(\theta)$} & $P_{k'}(\theta)$ & $P_{k'}(\theta')$ & \multirow{-2}{*}{\textbf{Ratio}} & \multirow{-2}{*}{$\varepsilon_{\mathrm{emp}}$} \\
  \midrule
  MLP & $0.998$ & $3.9 \times 10^{-7}$ & $2.4 \times 10^{-5}$ & $61\times$ & 4.11 \\
  VQC & $0.960$ & $7.6 \times 10^{-3}$ & $2.3 \times 10^{-2}$ & $3.0\times$ & \textbf{1.11} \\
  \bottomrule
  \end{tabular}
  \begin{tablenotes}
  \item[] {\footnotesize $P_{k^*}$: top-1 probability; $P_{k'}$: class maximizing $|\log(P_k/P_k')|$. Same input pair, $\|\theta - \theta'\|_2 = 0.05$.}
  \end{tablenotes}
  \end{threeparttable}
\end{table}

We examine specific adjacent input pairs to investigate this gap further and compare the prediction probabilities of both models. As shown in Table~\ref{tab:confidence_comparison}, the MLP exhibits severe overconfidence: its top-1 probability $P_{k^*}$ reaches $0.998$, forcing other class probabilities to near-zero (e.g., $P_{k'} < 10^{-6}$). When the input is slightly perturbed, these near-zero probabilities shift by orders of magnitude, causing the probability ratio to elevate ($61\times$) and yielding high $\varepsilon_{\mathrm{emp}}$. In contrast, the VQC maintains relatively moderate confidence ($P_{k^*} = 0.960$), with smoother probability distributions and less sensitivity to perturbations. This difference arises from the VQC's architectural constraint: its outputs are expectation values of Pauli-$Z$ measurements, inherently bounded in $[-1, +1]$, which acts as implicit regularization against overconfidence.

\begin{figure}[htbp]
    \centering
    \includegraphics[width=1\columnwidth]{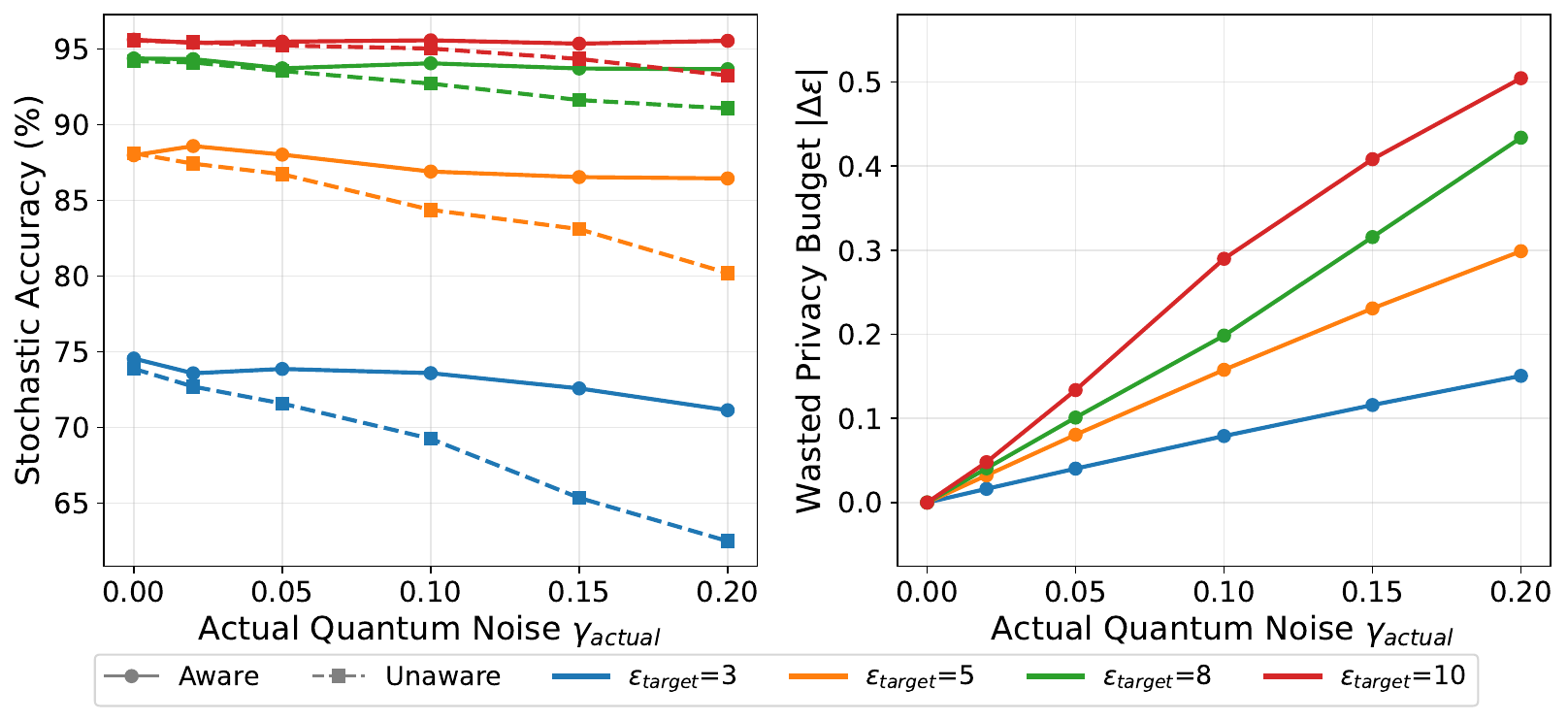}
    \caption{Noise-awareness benefit for privacy-utility tradeoff. Left: wasted privacy budget $|\Delta\varepsilon|$ due to unawareness. Right: stochastic accuracy comparison between noise-aware users (solid) and noise-unaware users (dashed) across privacy targets $\varepsilon \in \{3, 5, 8, 10\}$. Noise-aware users achieve up to 8.65\% higher accuracy at the same privacy level.}
    \label{fig:unawareness}
\end{figure}

\begin{table*}[htbp]
  \centering
  \caption{Performance Comparison at Different Privacy Budgets}
  \label{tab:eps_aligned_full}
  \resizebox{\linewidth}{!}{%
  \renewcommand{\arraystretch}{1.05}
  \begin{tabular}{c c c c c c c c c c c c c c c c c}
  \toprule
    & & \multicolumn{6}{c}{Generator Prediction Error (MW)} & \multicolumn{6}{c}{Power Flow Prediction Error (MW)} & & & \\
  \cmidrule(lr){3-8} \cmidrule(lr){9-14}
  \multirow{-2}{*}{Method} & \multirow{-2}{*}{$\varepsilon$} & $P_{g,1}$ & $P_{g,2}$ & $P_{g,3}$ & $P_{g,4}$ & $P_{g,5}$ & $P_{g,6}$ & $P_{l,10}$ &
  $P_{l,9}$ & $P_{l,8}$ & $P_{l,7}$ & $P_{l,6}$ & $P_{l,4}$ & \multirow{-2}{*}{\textbf{MAE}} & \multirow{-2}{*}{\textbf{Cost Gap}} & \multirow{-2}{*}{\textbf{Infeas. (\%)}} \\
  \midrule
  \rowcolor{gray!15}
  VQC & 0.9 & 0.533 & 0.751 & 0.201 & 0.083 & 0.280 & 0.594 & 0.721 & 0.721 & 0.721 & 0.552 & 0.552 & 0.552 & \textbf{0.522} & \textbf{4.96\%} & \textbf{23.6\%} \\
  MLP & 1.0 & 2.528 & 2.758 & 0.459 & 0.436 & 1.657 & 2.530 & 3.517 & 3.517 & 3.517 & 2.753 & 2.753 & 2.753 & 2.432 & 16.72\% & 67.3\% \\
  \midrule
  \rowcolor{gray!15}
  VQC & 2.8 & 0.030 & 0.038 & $<$0.01 & $<$0.01 & 0.019 & 0.041 & 0.054 & 0.054 & 0.054 & 0.042 & 0.042 & 0.042 & \textbf{0.036} & \textbf{0.21\%} & \textbf{4.6\%} \\
  MLP & 3.3 & 2.386 & 2.608 & 0.433 & 0.402 & 1.586 & 2.398 & 3.361 & 3.361 & 3.361 & 2.617 & 2.617 & 2.617 & 2.312 & 16.07\% & 66.1\% \\
  \midrule
  \rowcolor{gray!15}
  VQC & 5.0 & 0.014 & 0.018 & $<$0.01 & $<$0.01 & $<$0.01 & 0.023 & 0.030 & 0.030 & 0.030 & 0.023 & 0.023 & 0.023 & \textbf{0.019} & \textbf{0.09\%} & \textbf{2.2\%} \\
  MLP & 5.4 & 1.975 & 2.313 & 0.497 & 0.301 & 1.159 & 1.962 & 2.800 & 2.800 & 2.800 & 2.210 & 2.210 & 2.210 & 1.936 & 13.47\% & 63.9\% \\
  \midrule
  \rowcolor{gray!15}
  VQC & 9.7 & 0.012 & 0.017 & $<$0.01 & $<$0.01 & $<$0.01 & 0.021 & 0.027 & 0.027 & 0.027 & 0.021 & 0.021 & 0.021 & \textbf{0.017} & \textbf{0.08\%} & \textbf{1.8\%} \\
  MLP & 10.3 & 0.883 & 1.036 & 0.163 & 0.227 & 0.602 & 1.003 & 1.519 & 1.519 & 1.519 & 1.112 & 1.112 & 1.112 & 0.984 & 5.52\% & 45.1\% \\
  \bottomrule
  \end{tabular}%
  }
\end{table*}


To further explore the practical value of noise characterization, we examine the wasted privacy budget when users are unaware of the intrinsic quantum noise. Given a privacy target $\varepsilon_{\mathrm{target}}$ and actual noise level $\gamma_{\mathrm{actual}}$, we compute the required inverse temperature $\beta$ by inverting the bound in Theorem~\ref{thm:dp}: aware users substitute the true $\gamma_{\mathrm{actual}}$, while unaware users conservatively assume $\gamma=0$. Since intrinsic noise already contributes to privacy protection via the $(1-\gamma)$ factor, aware users need less additional uncertainty from softmax sampling and can use a larger $\beta$.

Fig.~\ref{fig:unawareness} quantifies the resulting inefficiency. The left figure shows the wasted privacy budget: unaware users unknowingly achieve $\varepsilon_{\mathrm{actual}} = \varepsilon_{\mathrm{target}}(1-\gamma) < \varepsilon_{\mathrm{target}}$, providing stronger protection than required. The wasted privacy budget grows nearly linearly with $\gamma_{\text{actual}}$, reaching $|\Delta\varepsilon| = 0.5$ at $\varepsilon_{\mathrm{target}}=10$ and $\gamma=0.2$. The right figure shows the corresponding accuracy loss: aware users (solid) consistently outperform unaware users (dashed), with the gap reaching nearly 10\% at $\varepsilon_{\mathrm{target}}=3$ and $\gamma=0.2$.

\subsection{Analysis of Operational Feasibility}

Practical deployment requires that the randomized region selection does not severely compromise operational feasibility or economic performance. In this section, we evaluate three metrics: prediction error for key decision variables (generators and line flows), constraint feasibility under the affine dispatch policy, and the economic performance in terms of cost gap.

We compute the affine solution map $x = F_{\tilde{k}} \theta + f_{\tilde{k}}$ for the predicted region $\tilde{k}$ and check whether all constraints $Wx \le S + T\theta$ are satisfied. The maximum violation threshold is $10^{-4}$. For cost evaluation, misclassified samples produce invalid cost estimates since each affine solution map is only optimal within its own region. We therefore project infeasible solutions onto the feasible set before computing the cost gap relative to the optimal solution.

Table~\ref{tab:eps_aligned_full} compares the VQC and MLP at matched privacy budgets. For each target $\varepsilon$, we select the $(\gamma, \beta)$ configuration that achieves the closest empirical $\varepsilon_{95}$ and evaluate the resulting dispatch quality. From the results, we can see that at comparable privacy levels, the VQC consistently achieves lower prediction errors, smaller cost gaps, and lower infeasibility rates. Specifically, at $\varepsilon \approx 10$, the VQC achieves MAE of 0.017~MW and infeasibility rate of 1.8\%, while the MLP yields MAE of 0.984~MW and infeasibility of 45.1\%. This advantage arises because the VQC can operate with a larger $\beta$ while meeting the same privacy target, resulting in more accurate region predictions.

\begin{figure}[htbp]
    \centering
    \includegraphics[width=\columnwidth]{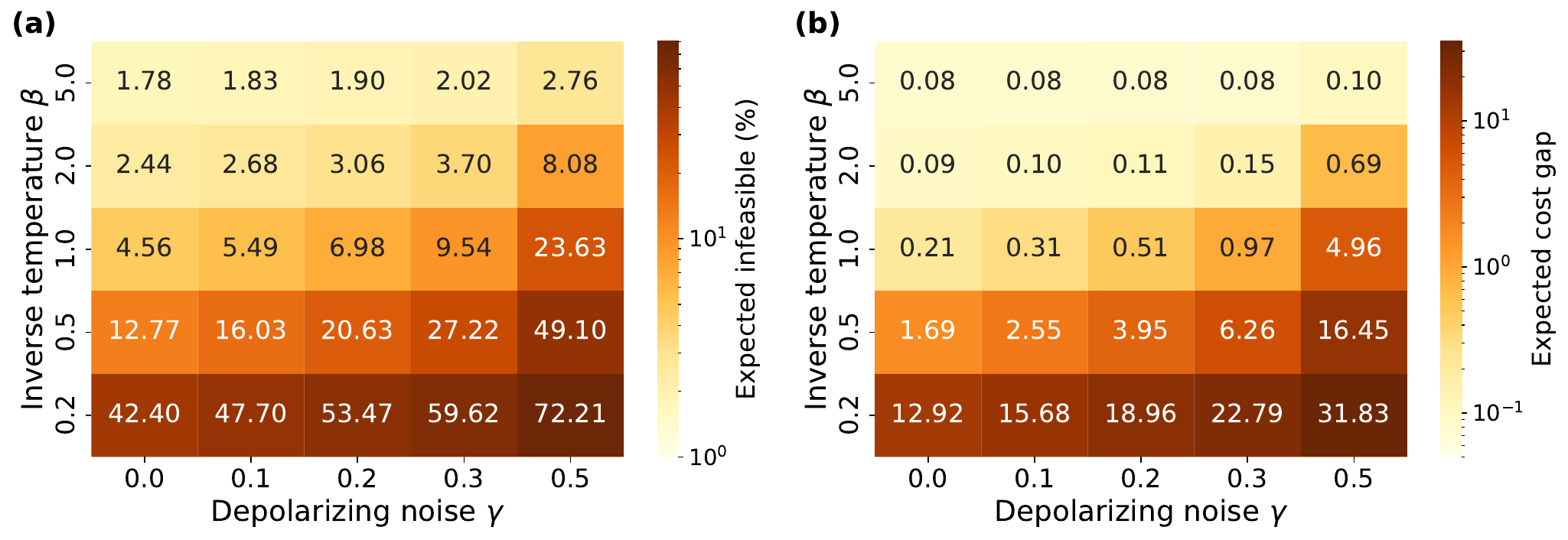}
    \caption{Infeasibility rate (\%) and cost gap across different noise levels $\gamma$ and inverse temperatures $\beta$. At high $\beta$ values, infeasibility remains stable regardless of noise level.}
    \label{fig:heatmap_infeasible}
\end{figure}

Fig.~\ref{fig:heatmap_infeasible} provides a detailed view of how noise level $\gamma$ and inverse temperature $\beta$ jointly affect operational metrics. The left figure shows expected infeasibility rates across the $(\gamma, \beta)$ grid. At typical operating points with $\beta \ge 2$, infeasibility remains below 5\% even at high noise levels such as $\gamma = 0.3$. The worst-case infeasibility of 72.2\% occurs at $(\gamma = 0.5, \beta = 0.2)$, corresponding to extremely strong privacy settings rarely needed in practice. The right figure shows expected cost gaps, which follow a similar pattern. Together, Fig.~\ref{fig:heatmap_infeasible} reveals that for large $\beta$, both metrics become nearly independent of noise level $\gamma$. At high $\beta$, softmax approximates argmax selection, and since depolarizing noise contracts all logits uniformly without changing the ranking, the selected region remains stable. This allows operators to add quantum noise for privacy without sacrificing operational reliability.

\subsection{Qubit Efficiency and Computational Speedup}

Practical deployment on near-term quantum hardware faces two constraints: limited qubit counts and computational efficiency for probabilistic analysis. This subsection analyzes how Algorithm 1 addresses both challenges. We first compare the qubit requirements between direct quantum optimization and our classification approach, then evaluate the computational speedup for POPF scenario evaluation.

\begin{table}[htbp]
    \centering
    \begin{threeparttable}
    \caption{Qubit Budget for Direct Quantum OPF vs. Region Classification}
    \label{tab:qubit_budget}
    \begin{tabular}{cc|cc|cc}
        \toprule
        \multicolumn{2}{c|}{Precision (bits)} & \multicolumn{2}{c|}{Qubits required} & \multicolumn{2}{c}{Total} \\
        \cmidrule(lr){1-2} \cmidrule(lr){3-4} \cmidrule(lr){5-6}
        Variable ($b$) & Slack ($Y$) & Variables & Slack & Direct QUBO & Ours \\
        \midrule
        \multirow{3}{*}{4} & 2 & 168 & 428  & 596  & \multirow{3}{*}{5} \\
                           & 3 & 168 & 642  & 810  & \\
                           & 4 & 168 & 856  & 1024 & \\
        \midrule
        \multirow{3}{*}{6} & 2 & 252 & 428  & 680  & \multirow{3}{*}{5} \\
                           & 3 & 252 & 642  & 894  & \\
                           & 4 & 252 & 856  & 1108 & \\
        \midrule
        \multirow{3}{*}{8} & 3 & 336 & 642  & 978  & \multirow{3}{*}{5} \\
                           & 4 & 336 & 856  & 1192 & \\
                           & 5 & 336 & 1070 & 1406 & \\
        \bottomrule
    \end{tabular}
    \begin{tablenotes}
    \item[] {\footnotesize Variables = $b \times 42$; Slack = $Y \times 214$ constraints. $Y$ depends on constraint range, typically $Y \geq \lceil b/2 \rceil$~\cite{morstyn2022annealing}.}
    \end{tablenotes}
    \end{threeparttable}
\end{table}

Direct quantum optimization algorithms such as quantum approximate optimization algorithm (QAOA) and quantum annealing solve problems formulated as quadratic unconstrained binary optimization (QUBO), which requires encoding continuous decision variables into binary representations~\cite{morstynOpportunitiesQuantumComputing2024}. For an OPF problem with $n$ continuous variables and $b$-bit precision, the qubit count scales as $\mathcal{O}(nb)$ for variable encoding, plus additional slack qubits for inequality constraints~\cite{morstyn2022annealing}. Table~\ref{tab:qubit_budget} compares these requirements for the 69-bus test system: decision variables require $b \times 42$ qubits, while the 214 inequality constraints require slack variables encoded with $Y$ bits each.

Beyond qubit feasibility, POPF requires evaluating thousands of uncertainty scenarios, making per-instance computation time critical. Table~\ref{tab:runtime} compares the per-instance runtime of different approaches on the 69-bus test system. The OPF solver (YALMIP with GUROBI) requires 13,833~$\mu$s per instance, measured on a standard desktop CPU.

The proposed approaches achieve substantial speedup by replacing optimization with classification plus affine evaluation. Constraint checking, which identifies the active region by testing all polyhedral boundaries, achieves 21.1~$\mu$s (656$\times$ speedup). The MLP classifier with affine solution map achieves 32.9~$\mu$s (420$\times$ speedup), dominated by neural network inference time.

For the VQC, we estimate runtime following the gate-model assumptions used in prior quantum power systems studies~\cite{zhuang2025quantum,zhou2024carbon}. The runtime is modeled as
\begin{equation}
    T_{\mathrm{VQC}} = T_P + T_G \times D + T_M,
\end{equation}
where $T_P + T_M = 1$~$\mu$s accounts for state preparation and measurement, $T_G = 10$~ns is the gate execution time, and $D$ is the circuit depth. For our 5-qubit, 6-layer circuit with depth $D = 1 + 6 \times (1 + 5) = 37$, this yields $T_{\mathrm{VQC}} = 1.37~\mu$s. Including affine solution map evaluation ($\sim$0.03~$\mu$s), the total is approximately 1.4~$\mu$s, corresponding to a speedup of 9,975$\times$. While this estimate assumes ideal gate-model hardware, it indicates the potential for accelerating large-scale POPF analysis with future quantum devices.

\begin{table}[htbp]
    \centering
    \caption{Runtime Comparison for Online OPF Evaluation}
    \label{tab:runtime}
    \begin{tabular}{ccc}
        \toprule
        Method                     & Runtime ($\mu$s) & Speedup       \\
        \midrule
        OPF Solver (YALMIP)        & 13,833.4         & 1$\times$     \\
        Constraint Check + Affine  & 21.1             & 656$\times$   \\
        MLP + Affine               & 32.9             & 420$\times$   \\
        VQC + Affine               & 1.4              & 9,975$\times$ \\
        \bottomrule
    \end{tabular}
\end{table}

\section{Conclusion}
\label{sec:conclusion}
This paper proposes a quantum-enabled probabilistic optimal power flow framework with built-in privacy guarantees and qubit efficiency. By reformulating POPF via multi-parametric programming, this framework encodes only the few uncertainty parameters rather than high-dimensional physical variables, reducing qubit requirements by two orders of magnitude. A variational quantum classifier identifies the active critical region for each scenario, with an estimated three orders of magnitude speedup in online evaluation. We prove that depolarizing noise and softmax sampling provide differential privacy guarantees as closed-form functions of measurable parameters. Case studies show the VQC achieves $2.1\times$ lower privacy budget
than classical neural networks, attributable to the bounded output structure of quantum measurements. At matched privacy levels, the VQC maintains lower infeasibility and prediction error. Future work includes validation on real quantum hardware.

\bibliographystyle{IEEEtran}
\bibliography{mybib}

\appendices
\setcounter{equation}{0}
\renewcommand{\theequation}{A.\arabic{equation}}
\section{Proof of Theorem~\ref{thm:dp}: Differential Privacy Guarantee}
\label{appendix:thm-dp}

The proof first bounds the trace distance between noisy quantum states under adjacent inputs, then propagates this bound through the measurements to the score vector, and finally applies the softmax sensitivity analysis to obtain the privacy guarantee. 

Recall that the trace distance $D_{\mathrm{tr}}(\rho, \rho') := \frac{1}{2}\|\rho - \rho'\|_1$ measures distinguishability between quantum states. For any matrix~$A$, $\|A\|_\infty$ denotes its largest singular value and $\|A\|_1$ the sum of its singular values (trace norm).

Fix adjacent parameters $\theta \sim \theta'$. Let $\rho_{\mathrm{enc}}(\theta) := \ket{\psi_{\mathrm{enc}}(\theta)}\!\bra{\psi_{\mathrm{enc}}(\theta)}$ denote the encoded input state and $\rho_{\mathrm{out}}(\theta) := U(\phi_q)\,\rho_{\mathrm{enc}}(\theta)\,U^\dagger(\phi_q)$ the circuit output state, and likewise for~$\theta'$. Since applying the same unitary gate to both states does not change their trace distance, we have
\[
D_{\mathrm{tr}}(\rho_{\mathrm{out}}(\theta),\,\rho_{\mathrm{out}}(\theta'))
=
D_{\mathrm{tr}}(\rho_{\mathrm{enc}}(\theta),\,\rho_{\mathrm{enc}}(\theta')).
\]
That is, the trainable circuit does not amplify the distinguishability introduced by the input encoding. By Assumption~\ref{ass:lipschitz},
\begin{equation}
D_{\mathrm{tr}}(\rho_{\mathrm{enc}}(\theta),\rho_{\mathrm{enc}}(\theta'))
\le L_{\mathrm{enc}}\|\theta-\theta'\|_2
\le L_{\mathrm{enc}}\Delta_\theta.
\end{equation}
Therefore,
\begin{equation}
D_{\mathrm{tr}}(\rho_{\mathrm{out}}(\theta),\rho_{\mathrm{out}}(\theta'))
\le L_{\mathrm{enc}}\Delta_\theta.
\end{equation}

Under the depolarizing channel,
\begin{equation}
    \mathcal{E}_\gamma(\rho) = (1-\gamma)\rho + \gamma \frac{I}{D},
\end{equation}
where $D = 2^{n_q}$ is the Hilbert-space dimension. We have,
\begin{equation}
\mathcal{E}_\gamma(\rho_{\mathrm{out}}(\theta))
-
\mathcal{E}_\gamma(\rho_{\mathrm{out}}(\theta'))
=
(1-\gamma)\bigl(\rho_{\mathrm{out}}(\theta)-\rho_{\mathrm{out}}(\theta')\bigr).
\end{equation}
Substituting this identity into the definition of trace distance ($D_{\mathrm{tr}}(\rho, \rho') := \frac{1}{2}\|\rho - \rho'\|_1$) gives
\begin{equation}
\begin{aligned}
    &\quad D_{\mathrm{tr}}\bigl(
        \mathcal{E}_\gamma(\rho_{\mathrm{out}}(\theta)),
        \mathcal{E}_\gamma(\rho_{\mathrm{out}}(\theta'))
    \bigr) \\
    &=
    \frac{1}{2}
    \left\|
        \mathcal{E}_\gamma(\rho_{\mathrm{out}}(\theta))
        -
        \mathcal{E}_\gamma(\rho_{\mathrm{out}}(\theta'))
    \right\|_1 \\
    &=
    \frac{1-\gamma}{2}
    \left\|
        \rho_{\mathrm{out}}(\theta)-\rho_{\mathrm{out}}(\theta')
    \right\|_1 \\
    &=
    (1-\gamma)
    D_{\mathrm{tr}}\bigl(
        \rho_{\mathrm{out}}(\theta),
        \rho_{\mathrm{out}}(\theta')
    \bigr) \\
    &\le (1-\gamma)L_{\mathrm{enc}}\Delta_\theta.
\end{aligned}
\end{equation}

Each feature coordinate is defined as the expectation value of the Pauli-$Z$ observable on qubit-$j$. Specifically,
\begin{equation}
    h_j(\theta)
    =
    \mathrm{Tr}\!\left(
        \sigma_z^{(j)}\,\mathcal{E}_\gamma(\rho_{\mathrm{out}}(\theta))
    \right).
\end{equation}
Hence, the difference in the $j$-th feature under adjacent inputs is
\begin{equation}
    h_j(\theta)-h_j(\theta')
    =
    \mathrm{Tr}\!\left(
        \sigma_z^{(j)}
        \left[
            \mathcal{E}_\gamma(\rho_{\mathrm{out}}(\theta))
            -
            \mathcal{E}_\gamma(\rho_{\mathrm{out}}(\theta'))
        \right]
    \right).
\end{equation}
By the inequality $|\mathrm{Tr}(AB)| \le \|A\|_\infty \|B\|_1$,
\begin{equation}
\begin{aligned}
    |h_j(\theta)-h_j(\theta')|
    &\le
    \|\sigma_z^{(j)}\|_\infty
    \left\|
        \mathcal{E}_\gamma(\rho_{\mathrm{out}}(\theta))
        -
        \mathcal{E}_\gamma(\rho_{\mathrm{out}}(\theta'))
    \right\|_1 \\
    &=
    2\|\sigma_z^{(j)}\|_\infty
    D_{\mathrm{tr}}\bigl(
        \mathcal{E}_\gamma(\rho_{\mathrm{out}}(\theta)),
        \mathcal{E}_\gamma(\rho_{\mathrm{out}}(\theta'))
    \bigr).
\end{aligned}
\end{equation}

Since the Pauli-$Z$ observable has eigenvalues $\pm 1$, we have $\|\sigma_z^{(j)}\|_\infty = 1$. Therefore,
\begin{equation}
    |h_j(\theta)-h_j(\theta')|
    \le
    2(1-\gamma)L_{\mathrm{enc}}\Delta_\theta.
\end{equation}
Since this holds for every coordinate $j$, it follows that
\begin{equation}
    \|h(\theta)-h(\theta')\|_\infty
    \le
    2(1-\gamma)L_{\mathrm{enc}}\Delta_\theta.
\end{equation}

For each critical region $k$, the logit is
\begin{equation}
    s_k(\theta)=w_k^\top h(\theta)+b_k,
\end{equation}
where $w_k^\top$ is the $k$-th row of $W_{cl}$. The bias term cancels in the difference, so
\begin{equation}
\begin{aligned}
    |s_k(\theta)-s_k(\theta')|
    &=
    \bigl|w_k^\top(h(\theta)-h(\theta'))\bigr| \\
    &\le
    \|w_k\|_1\,\|h(\theta)-h(\theta')\|_\infty \\
    &\le
    2(1-\gamma)L_{\mathrm{enc}}\Delta_\theta\,\|w_k\|_1.
\end{aligned}
\end{equation}
Define the logit sensitivity
\begin{equation}
    \Delta_s := \max_k |s_k(\theta)-s_k(\theta')|.
\end{equation}
Then
\begin{equation}
    \Delta_s
    \le
    2(1-\gamma)L_{\mathrm{enc}}\Delta_\theta\,\|W_{cl}\|_{\infty,1},
\end{equation}
where $\|W_{cl}\|_{\infty,1} := \max_k \|w_k\|_1$. 
Under softmax sampling, the probability of releasing region $k$ is
\begin{equation}
    \mathbb{P}[\tilde{k}=k\mid\theta]
    =
    \frac{\exp(\beta s_k(\theta))}{\sum_{j=1}^K \exp(\beta s_j(\theta))}.
\end{equation}
Therefore, for any region $k$, 
\begin{align}
    \frac{\mathbb{P}[\tilde{k} = k \mid \theta]}
         {\mathbb{P}[\tilde{k} = k \mid \theta']}
    &=
    \exp\!\bigl(\beta(s_k(\theta)-s_k(\theta'))\bigr)
    \cdot
    \frac{\sum_{j=1}^K \exp(\beta s_j(\theta'))}
         {\sum_{j=1}^K \exp(\beta s_j(\theta))}.
\end{align}
The first factor is at most $\exp(\beta\Delta_s)$ because
$|s_k(\theta)-s_k(\theta')| \le \Delta_s$ for every $k$. 
For the second factor, the same bound implies
$s_j(\theta') \le s_j(\theta)+\Delta_s$ for each $j$, and hence
\[
\exp(\beta s_j(\theta'))
\le
\exp(\beta\Delta_s)\exp(\beta s_j(\theta)).
\]
Summing over $j$ gives
\[
\sum_{j=1}^K \exp(\beta s_j(\theta'))
\le
\exp(\beta\Delta_s)\sum_{j=1}^K \exp(\beta s_j(\theta)).
\]
Therefore,
\begin{equation}
    \frac{\mathbb{P}[\tilde{k} = k \mid \theta]}
         {\mathbb{P}[\tilde{k} = k \mid \theta']}
    \le
    \exp(2\beta\Delta_s).
\end{equation}
Since the released output $\tilde{k}$ is discrete, this pointwise bound implies
$(\varepsilon_{\mathrm{reg}},0)$-differential privacy with
\[
\varepsilon_{\mathrm{reg}} = 2\beta\Delta_s.
\]
Substituting the bound on $\Delta_s$ gives
\begin{equation}
    \varepsilon_{\mathrm{reg}}
    \le
    4\beta(1-\gamma)L_{\mathrm{enc}}\Delta_\theta\,\|W_{cl}\|_{\infty,1},
\end{equation}
which is exactly~\eqref{eq:eps-bound}. \hfill$\square$

\ifincludeSecondAppendix
\setcounter{equation}{0}
\renewcommand{\theequation}{B.\arabic{equation}}
\section{Proof of Theorem~\ref{thm:tradeoff}: Privacy-Cost Trade-off}
\label{appendix:thm-tradeoff}

Let $x_k(\theta) = F_k \theta + f_k$ denote the affine solution for region $k$, and recall $\Delta J_k(\theta) = J(x_k(\theta); \theta) - J(x^*(\theta); \theta)$ with $\Delta J_{k^*}(\theta) = 0$.

\paragraph{Part (i): Cost Decomposition}
Let $\pi_k(\theta) := \mathbb{P}[\tilde{k} = k \mid \theta]$ and
$P_{\mathrm{err}}(\theta) := \mathbb{P}[\tilde{k} \ne k^*(\theta) \mid \theta] = 1 - \pi_{k^*}(\theta)$.
Since the realized cost gap is $\Delta J = \Delta J_{\tilde{k}}(\theta)$, the law of total expectation gives
\begin{equation}
    \mathbb{E}[\Delta J \mid \theta]
    = \sum_{k=1}^K \pi_k(\theta)\,\Delta J_k(\theta).
\end{equation}
Since $\Delta J_{k^*}(\theta) = 0$ and $\Delta J_k(\theta) \le \Delta J_{\max}(\theta)$ for all $k \ne k^*$, we obtain
\begin{equation}
\begin{aligned}
    \sum_{k=1}^K \pi_k(\theta)\,\Delta J_k(\theta)
    &= \sum_{k \ne k^*} \pi_k(\theta)\,\Delta J_k(\theta) \\
    &\le \Delta J_{\max}(\theta)\sum_{k \ne k^*} \pi_k(\theta) \\
    &= \Delta J_{\max}(\theta)\,P_{\mathrm{err}}(\theta).
\end{aligned}
\end{equation}

\paragraph{Part (ii): Margin Bound}
The probability of selecting the correct region under softmax sampling can be rewritten as
\begin{align}
    \pi_{k^*}(\theta)
    &= \frac{\exp(\beta s_{k^*}(\theta))}{\sum_k \exp(\beta s_k(\theta))} \\
    &= \frac{1}{1 + \sum_{k \ne k^*} \exp\bigl(\beta (s_k(\theta) - s_{k^*}(\theta))\bigr)}.
\end{align}
By the definition of the margin,
\begin{equation}
    m(\theta) = s_{k^*}(\theta) - \max_{k \ne k^*} s_k(\theta),
\end{equation}
we have $s_k(\theta) - s_{k^*}(\theta) \le -m(\theta)$ for all $k \ne k^*$. Therefore,
\begin{equation}
    \sum_{k \ne k^*} \exp\bigl(\beta (s_k(\theta) - s_{k^*}(\theta))\bigr)
    \le (K-1)\exp(-\beta m(\theta)),
\end{equation}
and hence
\begin{equation}
    \pi_{k^*}(\theta) \ge \frac{1}{1 + (K-1)\exp(-\beta m(\theta))}.
\end{equation}
It follows that
\begin{equation}
\begin{aligned}
    P_{\mathrm{err}}(\theta)
    &= 1 - \pi_{k^*}(\theta) \\
    &\le \frac{(K-1)\exp(-\beta m(\theta))}{1 + (K-1)\exp(-\beta m(\theta))} \\
    &\le (K-1)\exp(-\beta m(\theta)).
\end{aligned}
\end{equation}

Combining parts (i) and (ii) yields
\begin{equation}
\begin{aligned}
    \mathbb{E}[\Delta J \mid \theta]
    &\le \Delta J_{\max}(\theta)\,P_{\mathrm{err}}(\theta) \\
    &\le \Delta J_{\max}(\theta)\cdot (K-1)\exp\bigl(-\beta m(\theta)\bigr).
\end{aligned}
\end{equation} \hfill $\square$

\fi

\end{document}